\documentclass{IEEEojcsys}

\usepackage[colorlinks,urlcolor=blue,linkcolor=blue,citecolor=blue]{hyperref}

\usepackage{color,array}

\usepackage{graphicx}

\jvol{00}
\jnum{XX}
\paper{1234567}
\pubyear{2021}
\receiveddate{XX September 2021}
\accepteddate{XX October 2021}
\publisheddate{XX November 2021}
\currentdate{XX November 2021}
\doiinfo{OJCSYS.2021.Doi Number}

\graphicspath{ {./img/} }

\newcommand{\intra}{{\rm intra}} %
\newcommand{\ctr}{{\rm ctr}} %

\newcommand{\EP}{{\frac{1}{\varepsilon}}} %
\newcommand{\EPt}{{\frac{t}{\varepsilon}}} %

\newcommand\PC[1]{^{(#1)}}

\newcommand\scalemath[2]{\scalebox{#1}{\mbox{\ensuremath{\displaystyle #2}}}}

\newcommand\SmallMatrixNBla[1]{{%
		\small\arraycolsep=0.9\arraycolsep\ensuremath{\begin{matrix}#1\end{matrix}}}}
	
\newcommand{\blkdiag}{{\rm blkdiag}} %

\usepackage{mathtools}
\usepackage{bm}
\usepackage{bbm}

\usepackage{graphicx}
\usepackage{amssymb}

\usepackage{algorithm}
\usepackage{algpseudocode}

\newcommand{\BE}{\mathbb{E}} 
\newcommand{\BS}{\mathbb{S}} 
\newcommand{\BG}{\mathbb{G}} 
\newcommand{\R}{\mathbb{R}} 
\newcommand{\BONE}{\mathbbm{1}} 
\newcommand{\BT}{\mathbb{T}} 
\newcommand{\BV}{\mathbb{V}} 
\newcommand{\BB}{\mathbb{B}} 

\newcommand{\CP}{\mathcal{P}} %
\newcommand{\CR}{\mathcal{R}} %
\newcommand{\CC}{\mathcal{C}} %
\newcommand{\CG}{\mathcal{G}} %
\newcommand{\CE}{\mathcal{E}} %
\newcommand{\CM}{\mathcal{M}} %
\newcommand{\CV}{\mathcal{V}} %
\newcommand{\CD}{\mathcal{D}} %

\newcommand{\sign}{{\rm sign}}
\newcommand{\inter}{{\rm inter}}

\newcommand{\ctl}{{\rm ctl}}

\newcommand{\diag}{{\rm diag}} %

\usepackage{amsthm}

\newtheorem{theorem}{Theorem}

\newtheorem{lemma}{Lemma}
\newtheorem{corollary}{Corollary}

\theoremstyle{definition}
\newtheorem{remark}{Remark}

\newtheorem{definition}{Definition}
\newtheorem{assumption}{Assumption}

\newcommand*{\QE}{\hfill\ensuremath{\blacksquare}}	

\newcommand*{\QEDA}{\null \hfill\ensuremath{\triangle}}

\begin{document}

\sptitle{Article Category}

\title{Vibrational Stabilization of Cluster Synchronization in Oscillator Networks} 

\editor{This paper was recommended by Associate Editor F. A. Author.}

\author{Yuzhen Qin\affilmark{1} (Member, IEEE)}

\author{Alberto Maria Nobili\affilmark{2}} 

\author{Danielle S. Bassett\affilmark{3}  (Member, IEEE)}

\author{Fabio Pasqualetti\affilmark{1} (Member, IEEE)}

\affil{ Department of Mechanical Engineering, University of
	California, Riverside, CA, USA} 
\affil{ Perceptual Robotics Laboratory at the IIM Institute, Department of Excellence in Robotics and A.I., Scuola Superiore Sant’Anna, Pisa, Italy} 

\affil{ Department of
	Bioengineering, Department of Electrical \& Systems
	Engineering,  Department of Physics \& Astronomy, the
	Department of Psychiatry, and  Department of Neurology,
	University of Pennsylvania, and the Santa Fe Institute}

\corresp{CORRESPONDING AUTHOR: Yuzhen Qin (e-mail: \href{mailto:yuzhenq@ucr.edu}{yuzhenq@ucr.edu})}
\authornote{This research was funded in part by Aligning Science Across	Parkinson’s (ASAP-020616) and in part by NSF (NCS-FO-1926829) and ARO (W911NF1910360).}

\markboth{PREPARATION OF PAPERS FOR IEEE OPEN JOURNAL OF CONTROL SYSTEMS}{F. A. AUTHOR {\itshape ET AL}.}

\begin{abstract}
	Cluster synchronization is of paramount importance for the normal functioning of numerous technological and natural systems. Deviations from normal cluster synchronization patterns are closely associated with various malfunctions, such as neurological disorders in the brain. Therefore, it is crucial to restore normal system functions by stabilizing the appropriate cluster synchronization patterns. Most existing studies focus on designing controllers based on state measurements to achieve system stabilization. However, in many real-world scenarios, measuring system states, such as neuronal activity in the brain, poses significant challenges, rendering the stabilization of such systems difficult. To overcome this challenge, in this paper, we employ an open-loop control strategy, \textit{vibrational control}, which does not requires any state measurements. We establish some sufficient conditions under which vibrational inputs stabilize cluster synchronization. Further, we provide a tractable approach to design vibrational control. Finally, numerical experiments are conducted to demonstrate our theoretical findings.
\end{abstract}

\begin{IEEEkeywords}
	Vibrational Control, Cluster Synchronization, Oscillator Networks
\end{IEEEkeywords}

\maketitle

\section{Introduction}
Cluster synchronization describes the phenomenon in which units within a network exhibit synchronized behavior, forming distinct clusters. This intriguing phenomenon is widely observed in many natural and engineering systems. For instance, it manifests as correlated neural activity in the brain. Different patterns of cluster synchronization in the brain play a fundamental role in various functions, including neuronal communication, memory formation and retrieval, and motor function \cite{HG-PAA-DG-AA-KA:19,FJ-AN:2011}.  

However, many brain disorders, such as Parkinson's disease \cite{Hammond2007} and epilepsy \cite{JP-DCM-JJGR-SCA:2013}, are characterized by aberrant synchrony patterns of brain activity. It becomes crucial to be able to stabilizing normal patterns of cluster synchronization. Most existing studies rely on the assumption that the states of systems can be measured to design feedback controllers to stabilize them. Unfortunately, many real-world systems often exhibit complex nonlinear dynamics over large network structures with states that are difficult to observe or measure in real time (such as neuronal activity in the brain), thus preventing the use of
sophisticated feedback techniques.

In this paper we leverage \emph{vibrational control} strategies to
ensure stability of network systems. Vibrational control is a powerful
strategy applicable in various domains to stabilize the dynamics of
complex systems without the need for direct state measurements. Unlike
traditional feedback-based control methods that rely on directly
measuring the system's states or outputs, vibrational control
leverages the inherent dynamics of the system to induce desired
stability and performance. In particular, vibrational control uses
pre-designed high-frequency signals, injected at specific locations
and times. As we show below, these signals can effectively change the
system dynamics and suppress unstable dynamics. Successful
applications of vibrational control to centralized systems are
numerous, including to inverted pendulums, chemical reactors, and
under-actuated robots \cite{REB-JB-SMM:86b,BS-BTZ:97,CX-TY-MI:2018}.
This paper takes the first steps to develop a theory of vibrational
control for network systems. Interestingly, in addition to its
technological value, the theory of vibrational control for network
systems may also help explain the success of deep brain stimulation
methods, as this technique also relies on the injection of
high-frequency electrical pulses to regulate brain processes and
restore healthy functions \cite{Krauss2021}. 


We focus on the stabilization of networks of heterogeneous Kuramoto
oscillators, and in particular on the stabilization of the dynamics
around a desired cluster synchronization manifold
\cite{TM-GB-DSB-FP:18}. We remark that (i) networks of Kuramoto
oscillators have been used successfully to model different phenomena
in diverse domains ranging from power engineering to biology and
neuroscience, thus making our results of broad applicability, and (ii)
cluster synchronization has been used as a proxy to model and regulate
the emergence of functional activation patterns in the brain
\cite{TM-GB-DSB-FP:19b,TM-GB-DSB-FP:22}, thus making our results of
timely relevance to these problems. 

\noindent
\textbf{Related work.}
Cluster synchronization has garnered significant attention recently as researchers seek to understand its underlying mechanisms and control strategies. Investigations into the field have revealed intriguing connections between cluster synchronization and network symmetries \cite{Pecora2014, YSC-TN-AEM:17, YQ-MC-BDOA-DSB-FP:20,EJ-SA-SJ-CJP-DRM:20} as well as equitable partitions \cite{Schaub2016}. Furthermore, stability conditions for cluster synchronization have been established in networks featuring dyadic connections \cite{TM-GB-DSB-FP:18, QY-KY-PO-CM:21, FP-SA-MT:21, FP-SA-MT:2020, KR-IH:2021} and hyper connections \cite{SA-DRM:2021a, SA-DRM:2021b}.
To control cluster synchronization, researchers have proposed diverse strategies, such as pinning control \cite{WW-WZ-TC:09} and interventions that involve manipulating network connections or the dynamics of individual nodes \cite{GLV-FM-LV:2018, DF-GR-MDB:17, TM-GB-DSB-FP:21-arxiv}. 
In contrast, our approach focuses on vibrational control, which offers a more realistic strategy in many real-world systems, e.g., for regulating neural activity as it resembles deep brain stimulation \cite{Krauss2021}. To the best of our knowledge, our work is among the first ones to utilize vibrations to regulate network systems.


\noindent
\textbf{Paper contributions.} The main contributions of this paper are
as follows. First, we formalize the problem of vibrational control for Kuramoto-oscillator networks, with the aim to stabilize patterns of cluster synchronization. By employing an averaging argument, we demonstrate that introducing vibrational inputs into the network effectively alters the system dynamics on average. We establish sufficient conditions for the effectiveness of vibrational control in stabilizing cluster synchronization within Kuramoto-oscillator networks. Through the analysis of linearized systems, we gain deep insights into the underlying mechanisms of vibrational control, revealing its ability to enhance the robustness of synchronization within clusters. Second, we establish a connection between the design of vibrational control in Kuramoto-oscillator networks and linear network systems. We show that vibrational control can effectively modify the weights of edges in linear network systems, thereby shaping their robustness. We develop graph-theoretical conditions for determining which edges can be modified through vibration. Additionally, we propose a systematic approach to designing vibrational control that targets the modifiable edges, aiming to enhance the overall robustness of the network system. Building on these findings, we apply the results to homogeneous Kuramoto-oscillator networks, presenting a tractable approach to designing vibrational control for improving the robustness of full synchronization. Furthermore, we extend the application of these results to heterogeneous Kuramoto-oscillator networks, deriving precise forms and placements of vibrational inputs to stabilize cluster synchronization. Finally, we conduct a numerical experiment to demonstrate our method for designing vibrational control to stabilize cluster synchronization in Kuramoto-oscillator networks.


A preliminary version of this work appeared in \cite{YQ-DSB-FP:22a}. With respect to \cite{YQ-DSB-FP:22a}, this paper (i) obtains a tighter sufficient condition for vibrational inputs to stabilize cluster synchronization, and (ii) presents a much more comprehensive approach to design vibrations in Kuramoto-oscillator network by deriving and utilizing precise forms and placements of vibrational inputs in linear network systems.

\textbf{Notation.} Given a matrix $M\in\R^{m\times n}$, the matrix $A = \sign(M):=[a_{ij}]$ is defined in a way such that $a_{ij}=1$ if $m_{ij}>0$, $a_{ij}=-1$ if $m_{ij}<0$, and $a_{ij}=0$ if $m_{ij}=0$. Denote the unit circle by $\BS^1$, and a point on it is called a \textit{phase} since the point can be used to indicate the phase angle of an oscillator.  For any two phases
$\theta_1,\theta_2 \in \BS^1$, the geodesic distance between them is the minimum of the lengths of the counter-clockwise and clockwise arcs connecting them, which is denoted by $|\theta_1-\theta_2 |_\BS$. Given a matrix $A$, $A^\dagger$ denote its pseudo-inverse.

\section{Problem Formulation}
\subsection{Kuramoto-Oscillator Networks}

Consider a Kuramoto-oscillator network described by 
\begin{equation}\label{no_input}
	\dot \theta_i =\omega_i + \sum_{j=1}^{n} w_{ij} \sin(\theta_j- \theta_i), 
\end{equation}
where $\theta_i \in \BS^1$ is the $i$th oscillator's phase, $\omega_i\in \R$ is its natural frequency for $i=1,\dots,n$, and $w_{ij}$ is the coupling strength. This paper investigates a network configuration where connections between nodes are bidirectional but allow for asymmetry (i.e., $w_{ij}\neq w_{ji}$). This assumption represents a more flexible scenario compared to the commonly studied undirected networks in existing literature. We use a weighted directed graph  $\CG=(\CV,\CE,W)$, to describe the network structure,  where $\CV=\{1,2,\dots, n\}$, $\CE\subseteq \CV\times \CV$, and $W=[w_{ij}]_{n\times n}$ with $w_{ij}\ge 0$ is the weighted adjacency matrix. There is a directed edge from $j$ to $i$, i.e., $(j,i)\in \CE$, if $w_{ij}>0$; $(j,i)\notin \CE$, otherwise.


In this paper, we are interest in studying cluster synchronization in the Kuramoto network \eqref{no_input}. Let us first provide a formal definition of cluster synchronization.

\begin{definition}[\textit{Cluster Synchronization manifold}]\label{CS:manifold}
	For the the network $\CG=(\CV,\CE)$, consider the partition $\CP:=\{\CP_1,\CP_2,\dots,\CP_r\}$, where each $\CP_k$ is a subset of $\CV$ satisfying $\CP_k \cap  \CP_\ell =\emptyset$ for any $k\neq \ell$, and $\cup_{k=1}^r \CP_k=\CV$. The \textit{cluster synchronization manifold} on the partition $\CP$ is defined as
	\begin{equation*}
		\CM:= \{\theta \in \BT^n: \theta_i=\theta_j, \forall i,j \in \CP_k,  k= 1,\dots,r\}.\hspace{12pt}\QEDA
	\end{equation*}
\end{definition} 

The manifold $\CM$ is said to be invariant along the system \eqref{no_input} if starting from $\theta(0)\in \CM$, the solution to \eqref{no_input} satisfies $\theta(t)\in \CM$ for all $t\ge 0$. We make the following assumption to ensure that $\CM$ is invariant. Note that similar assumptions are made in \cite{TM-GB-DSB-FP:18,LT-CF-MI-DSB-FP:17} for undirected networks.

\begin{assumption}[\textit{Invariance}]\label{invariance}
	For $k=1,2,\dots, r$:  i) the natural frequencies satisfy $\omega_i=\omega_j$ for any $i,j \in \CP_k$; and ii) the coupling strengths satisfy that, for any $\ell\in\{1,2,\dots,r\}\backslash\{k\}$,  $\sum_{q\in \CP_\ell}(a_{iq}-a_{jq})=0$ for any $i,j \in \CP_k$. \QEDA
\end{assumption}

To ensure the cluster synchronization defined by $\CM$ to appear, the stability of $\CM$ is also required.  Given a manifold $\mathcal C \in \BT^n$, define a $\delta$-neighborhood of $\mathcal C$ by $U_{\delta}(\mathcal C) = \{\theta \in \BT^n:{\rm dist}(\theta,\mathcal C)<\delta\}$ with ${\rm dist}(\theta,\mathcal C)=\inf _{y\in \mathcal C}{\|\theta-y\|_\BS}$. The  exponential stability of the manifold $\CM$ is defined below. 

\begin{definition}
	The manifold $\CM\in \BT^n$ is said to be exponentially stable  along the system \eqref{no_input}  if there is $\delta>0$ such that for any initial phase $\theta(0) \in\BT^ n$ satisfying $\theta(0)\in U_{\delta}( \CM)$ it holds that for all $t\ge 0$, ${\rm dist}(\theta(t),\CM)=k\cdot{\rm dist}(\theta(0), \mathcal C)\cdot e^{-\lambda t}$ for some $k>0$ and $\lambda>0$.
\end{definition}

Sufficient conditions are constructed for the exponential stability of $\CM$ (e.g., see \cite{TM-GB-DSB-FP:18,QY-KY-PO-CM:21,Schaub2016}). However, variations in network parameters due to factors like aging or brain disorders can disrupt these conditions, resulting in the loss of stability in cluster synchronization. This paper focuses on the application of vibrational control, a control strategy reminiscent of deep brain stimulation (DBS) \cite{Krauss2021}, to restore the stability of desired cluster synchronization patterns. Next, we will present an introduction to vibrational control.

\subsection{Vibrational Control}
Following \cite{REB-JB-SMM:86a}, consider a nonlinear system
\begin{equation}\label{general}
	\dot x = f(x, a),
\end{equation} 
where $x\in\R^n$, and $a\in \R^{m}$ is the parameter of the system. For linear systems, $a$ can be rewritten into a matrix form $A\in \R^{n\times n}$, and then we have $\dot x =A x$.  Vibrational control introduces vibrations to the parameters of \eqref{general}, leading to
\begin{equation}\label{general:ctred}
	\dot x = f(x, a+v(t)).
\end{equation} 
The control vector $v(t)=[v_i]\in \mathbb{R}^n$ is often selected to have the following structure:
\begin{equation}\label{vib:gener_form}
	v_{i}(t)= \sum_{\ell=1}^{\infty} \alpha_{i}^{(\ell)} \sin(\ell \beta_{i}t+\phi_{i}^{(\ell)}).
\end{equation}
which is periodic, zero-mean, and high-frequency \cite{REB-JB-SMM:86b}. Vibrational control is an open-loop strategy; an appropriate configuration of vibrations can stabilize an unstable system \eqref{general} without any measurements of the states \cite{SMM:80,BS-BTZ:97,CX-TY-MI:2018}. 

\subsection{Vibrational Control in Kuramoto-Oscillator Networks}
In the Kuramoto-oscillator network described by \eqref{no_input}, the natural frequencies and coupling strengths are the parameters. This paper specifically focuses on injecting vibrations solely into the couplings. Specifically, the control exerts its influence on the system by
\begin{equation}\label{Kuramoto:controlled}
	\dot \theta_i =\omega_i + \sum_{j=1}^{n} \big(w_{ij}+v_{ij}(t)\big) \sin(\theta_j- \theta_i),
\end{equation}
where $v_{ij}(t)$ is the vibration introduced to the edge $(j,i)$. Particularly, we consider that each $v_{ij}(t)$ is simply sinusoidal, which naturally satisfies \eqref{vib:gener_form}, i.e.,
\begin{equation}\label{form:vib}
	v_{ij}(t)= \mu_{ij } \sin\Big(\frac{\beta_{ij}}{\varepsilon} t\Big),
\end{equation}
Here, $\mu_{ij} \in \mathbb{R}^n$ represents the amplitude of the vibrations, while $\varepsilon>0$ is a small constant shared by all the vibrations, indicating their high-frequency nature. Additionally, $\beta_{ij}>0$ allows for distinct frequencies of vibrations injected into different edges. Note that various types of vibrations can be utilized in practical applications. However, for the sake of analysis simplification, we consider sinusoidal vibrations.

In contrast to the general system \eqref{general}, where vibrations can be applied to any parameter in $a$, the introduction of vibration control in a network system is subject to the constraints imposed by the network structure. It is reasonable to assume that vibrations can only be introduced to edges that already exist in the network. Therefore, vibrational control satisfies
\begin{align}\label{constraints}
	&\forall i,j: \mu_{ij}=0 &\text{ if } w_{ij}=0.
\end{align}

\textbf{Objective:} The objective of this paper is to design vibrational control inputs satisfying both \eqref{form:vib} and \eqref{constraints} to stabilize the cluster synchronization represented by the manifold $\mathcal{M}$.

\section{Preliminary}
\subsection{Graph-theoretical Notations}
\begin{figure}[t]
	\centering
	\includegraphics[scale=1.6]{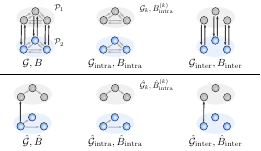}
	\caption{Summary of the main notations we use in this paper: different subgraphs and their corresponding incidence matrices. }
	\label{notation}
\end{figure}

We first define some graph-theoretical notations, which are further elucidated in a more intuitive manner in Figure~\ref{notation}.

For the directed graph $\CG=(\CV,\CE,C)$ that described the network in \eqref{no_input}, denote the oriented incidence matrix as $B=[b_{k\ell}]\in\R^{n\times 2m}$, where $\ell$ is the index of the edge $(i,j)\in \CE$; the elements in $B$ are defined in way where $b_{k\ell}=1$ if $k$ is the source of this edge, $b_{k\ell}=-1$ if $k$ is the sink of this edge, and $b_{k\ell}=0$ otherwise. 

For the partition $\CP:=\{\CP_1,\CP_2,\dots,\CP_r\}$ of $\CG$, define $\CG_k=(\CP_k,\CE_k)$ where $\CE_k:=\{(i,j)\in \CE: i,j\in \CP_k\}$. For each $k$, denote $n_k:=|\CP_k|$ as the number of nodes in $\CG_k$. We assume that each $\CG_k$ is connected and contains at least 2 nodes. Let $\CG_\intra=(\CV,\CE_\intra)=\cup_{k=1}^r \CG_k$ and we refer to it as the\textit{ intra-cluster sub-network}. Let $\CG_{\inter}=(\CV,\CE_\inter)$ be the \textit{inter-cluster sub-network}, where $\CE_\inter :=\CE\backslash \CE_\intra$. Let $B_\intra$ and $B_\inter$ be the oriented incidence matrices of $\CG_\intra$ and $\CG_\inter$, respectively.

Let $\hat \CG=(\CV,\hat \CE)$ be any spanning tree of $\CG$, and thus $|\hat \CE|=n-1$. Let $\hat B$ be the incidence matrix of $\hat \CG$.
For each $k$, Let $\hat \CG_{k}=(\CP_k,\hat \CE_k)$ with  $\hat \CE_k:=\hat \CE\cap \CE_k$. The intra-cluster subgraph of  the spanning tree $\hat \CG$ is $\hat \CG_\intra=(\CV,\hat \CE_\intra):=\hat \cup_{k=1}^r \hat \CG_k$;  the inter-cluster subgraph  of $\hat \CG$ is $\hat \CG_\inter=(\CV,\hat \CE_\inter)$ with $\hat \CE_\inter:=\hat \CE\backslash \hat \CE_\intra$. Denote $\hat B_\intra$ and $\hat B_\inter$ be the incidence matrices of  $\hat \CG_\intra$ and $\hat \CG_\inter$, respectively. 

Without loss of generality, we order the edges in the above incidence matrices in a way such that  
\begin{equation}\label{matrices}
	\begin{matrix*}[l]
			B=[B_\intra, B_\inter],&B_\intra=\blkdiag( B^{(1)}_{\intra},\dots, B^{(r)}_{\intra})\\
			\hat B=[\hat B_\intra, \hat B_\inter],&\hat B_\intra=\blkdiag(\hat B^{(1)}_{\intra},\dots,\hat B^{(r)}_{\intra}),
	\end{matrix*}
\end{equation}
where $B^{(k)}_{\intra}$ and $\hat B^{(1)}_{\intra}$ are the incidence matrices of the corresponding subgraph within individual clusters. 

\subsection{Incremental Dynamics}
To analyze the stability of the system \eqref{Kuramoto:controlled}, we define the incremental dynamics
\begin{align}\label{incremental}
	\dot \theta_k - \dot \theta_\ell =\omega_k -\omega_\ell &+ \sum_{j=1}^{n} \big(w_{ij}+v_{kj}(t)\big) \sin(\theta_j- \theta_k) \nonumber \\
	&- \sum_{j=1}^{n} \big(w_{\ell j}+v_{\ell j}(t)\big) \sin(\theta_j- \theta_\ell),
\end{align}
where each pair of $k$ and $\ell$ satisfies $(k,\ell)\in \hat \CE$. In other words, we only look at the dynamics of phase differences along the edges in the spanning tree. 

Let $x=\hat B_\intra \theta \in\R^{n-r}$ and $y=\hat B_\inter \theta \in \R^{r-1}$, which captures the intra- and inter-cluster phase differences. Then, one can derive that (see Appendix~\ref{derivatin} for the details)
\begin{subequations}\label{compact_form}
\begin{align}
	&\dot x=f_\intra(x) +f_\inter(x,y)+ h_\intra(V_\intra(t),x) \nonumber\\  
	&\hspace{3.5cm} +h_\inter\big(V_\inter(t),x,y\big),\label{compact_form:1}\\
	&\dot y=g(x,y)+h'_\intra(V_\intra(t),x) \nonumber\\  
	&\hspace{3.5cm} +h'_\inter\big(V_\inter(t),x,y\big)
\end{align}
\end{subequations}
where $x\in \R^{n-r}$, $y\in \R^{r-1}$, $V(t)=\diag ([v_{ij}]_{(i,j)\in \CE})\in \R^{2m}$ captures the vibrations introduced to the edges in the network, and the functions are given in \eqref{expr:function} in Appendix~\ref{derivatin},  describing the dynamics induced by the intra connections, inter connections, and vibrational control inputs, respectively. Here, we single out some important properties of these functions: i)  $f_\intra(0)=0$, ii) $f_\inter(0,y)=0$ for any $y$, iii) $h_\intra(V,0)=0$ iv) $h'_\intra(V,0)=0$ for any $V$.

Further, $x=0$ corresponds to the cluster synchronization described by $\CM$. However, to ensure $x=0$ is an equilibrium of \eqref{compact_form:1}, the vibrations need to satisfy a similar condition to Assumption~\ref{invariance}. To meet this requirement and also simplify analysis, we assume that vibrations are only injected to intra-cluster edges, i.e., $V_\inter(t)=0$. As a result, $h_\inter =0$ and $h'_\inter=0$. Then, the system \eqref{compact_form} becomes
\begin{subequations}\label{intra_only}
	\begin{align}
		&\dot x=f_\intra(x) +f_\inter(x,y)+  h_\intra(V_\intra(t),x),\label{intra_only:1}\\
		&\dot y=g(x,y)+ h'_\intra(V_\intra(t),x).
	\end{align}
\end{subequations}

Recall that vibrations usually need to be high-frequency. To streamline the analysis and simplify notation, we denote $f_\ctl := h_\intra$, $g_\ctl :=h'_\intra$. Also, for \eqref{form:vib}, we let $\mu_{ij}=\frac{u_{ij}}{\varepsilon}$ and $u_{ij}(t)=u_{ij}\sin(\beta_{ij}t)$. Let $U(t)=\diag ([u_{ij}(t)]_{(i,j)\in \CE})\in \R^{2m}$.  Then, the system \eqref{intra_only} can be rewritten into
\begin{subequations}\label{compact:main}
	\begin{align}
		&\dot x=f_\intra(x) +f_\inter(x,y)+ \EP f_\ctl\Big(U\big(\EPt \big),x \Big),\label{compact:main:1}\\
		&\dot y=g(x,y)+\EP g_\ctl(U\big(\EPt\big),x).
	\end{align}
\end{subequations}

Consequently, the manifold $\CM$ is exponentially stable along the system \eqref{Kuramoto:controlled} if $x=0$ is partially exponentially stable along \eqref{compact:main}. To stabilize the cluster synchronization manifold $\CM$, it suffices to design vibration control inputs to ensure the partial exponential stability of  $x=0$. 


\section{Vibrational Stabilization: General Results}\label{results:general}
Notice that $f_\inter(0,y)=0$ holds for any $y$, then the term $f_\inter(x,y)$ in \eqref{compact:main:1} can be viewed as a vanishing perturbation dependent of $y$ to the controlled nominal system
\begin{align}\label{nominal}
	\dot x=f_\intra(x)+\EP f_\ctr( U(\EPt),x).
\end{align}
From \eqref{expr:inter-pert}, this perturbation can be decomposed as 
\begin{align*}
	f_\inter=[(f^{(1)}_{\inter})^\top,\dots,(f^{(r)}_{\inter})^\top]^\top
\end{align*}
where $f^{(k)}_{\inter}=- (\hat B^{(k)}_{\intra})^\top \BB_\inter W_\inter \sin(R_2 x+R_3y)$ is the perturbation received by the $k$th cluster\footnote{As defined in Appendix~\ref{derivatin},  $\BB$ is obtained by replacing the positive elements in the oriented incidence matrix $B$ by $0$. Let $\BB=[\BB_\intra,\BB_\inter]$, and $\BB_\intra=[\BB_\intra\PC{1},\BB_\intra\PC{2},\dots,\BB_\intra\PC{r}]$. Also, $R_2$ and $R_2$ are given in Appendix~\ref{derivatin}.}.

Next, we show how a vibration control $f_\ctr(U,x)$ can stabilize $x=0$ in the presence of the perturbation $f_\inter(x,y)$. To this end, we linearize the system at $x=0$ and obtain
\begin{align}\label{linearized}
	\scalemath{1}{\dot x = (J + \EP P(\EPt))x+ N(y)x},
\end{align}
where
\begin{equation}\label{expres:JandP}
	\begin{aligned}
		&{J=\frac{\partial f_\intra}{\partial x}(0)= \blkdiag (J\PC{1},\dots,J\PC{r})},  \\
		& {P(t)=\frac{\partial f_\ctr}{\partial x}(0)=  \blkdiag (P\PC{1}(t),\dots,P\PC{r}(t))},\\
		&N(y) = \frac{\partial f_\inter}{\partial x}(0,y)=\blkdiag (N\PC{1}(t),\dots,N\PC{r}(t)).
	\end{aligned} 
\end{equation}
Here,  for any $k$,  ${J\PC{k}=-(\hat B^{(k)}_{\intra})^\top \BB_\intra^{(k)} W^{(k)}_\intra R_1}$ and 
\begin{align*}
	&P\PC{k}(t)=-	(\hat B_\intra^{(k)})^\top \BB^{(k)}_\intra U^{(k)}(t) R_1.\\
	&N\PC{k}(y) = -(\hat B\PC{k}_\intra)^\top \BB_\inter W_\inter (\BONE_{r-1}\otimes\sin(R_3y))\odot R_2.
\end{align*}

Observe that $P(t)$ is periodic and has the same period $T$ as $U(t)$.

\begin{lemma}[\textit{Connecting the stability of Systems} \eqref{compact_form} \textit{and} \eqref{linearized}] \label{lemma:linearized}
	If the equilibrium $x=0$ is exponentially stable uniformly in $y$ for the system \eqref{linearized}, it is also exponentially stable uniformly in $y$  for the system  \eqref{compact_form}. 
\end{lemma}

The proof can be found in Appendix~\ref{proof:linearized}. From this lemma, one can see that to stabilize the cluster synchronization manifold $\CM$, it suffices to configure the vibrational control such that $P(t)$ stabilizes  $x=0$ of the system \eqref{linearized}. 

Let $s={t}/{\varepsilon}$. The system \eqref{linearized} can be rewritten as 
\begin{align}\label{change:time:scale}
	\frac{d x}{d s} = (\varepsilon J +  P(s))x + \varepsilon N(y)x.
\end{align}

Next, we use averaging methods to analyze this system. However, the standard first-order averaging is not applicable here. Recall that $P(s)$ has zero mean. Then, applying the first-order averaging to \eqref{change:time:scale} just eliminates the $P(s)$ term and results in the uncontrolled system $\frac{d x}{d s} = \varepsilon J x+\varepsilon P( y)x$. 

To avoid that, we change the coordinates of  \eqref{change:time:scale} first before using averaging method. To do that, we introduce an auxiliary system  
\begin{align}\label{periodic}
	{\frac{d \hat x}{d s}=P(s) \hat x},
\end{align}
and let $\Phi(s,s_0)$ be its state transition matrix. Since $P$ is block-diagonal, it holds that
$
\Phi=\blkdiag(\Phi\PC{1},\dots,\Phi\PC{r}),
$
where $\Phi\PC{k}$ is the transition matrix of the subsystem in the $k$th cluster ${d \hat x_k}/{d s}=P\PC{k}(s) \hat x_k$. 

Consider the change of coordinates $z(s)=\Phi^{-1}(s,s_0) x(s)$. It follows from the system \eqref{change:time:scale} that
\begin{align}\label{coordinated}
	\frac{dz}{ds} = \varepsilon \Phi^{-1} J \Phi z+ \varepsilon \Phi^{-1} N(y) \Phi z.
\end{align}
Since $P(s)$ is $T$-periodic, $\Phi$ and $\Phi^{-1}$ are also $T$-periodic.
Then, we associate \eqref{coordinated} with a partially averaged system
\begin{align}\label{average}
	{\frac{dz}{ds} = \varepsilon (\bar J + \Phi^{-1} N(y) \Phi)z},
\end{align}
where
\begin{align}\label{J-bar}
	{\bar J= \frac{1}{T}\int_{s_0}^{s_0+T} \Phi^{-1} (s,s_0) J \Phi(s,s_0)ds}.
\end{align}

As both $J$ and $\Phi$ are block-diagonal, one can derive that $\bar J$ is also block-diagonal satisfying $\bar J=\blkdiag(\bar J\PC{1},\dots,\bar J\PC{r})$ with 
\begin{align*}
	\bar J \PC{k}= \frac{1}{T}\int_{s_0}^{s_0+T} \Big(\Phi\PC{k}(s,s_0)\Big)^{-1} J\PC{k} \Phi\PC{k}(s,s_0)ds.
\end{align*} 
Recall that $\Phi$ and $\Phi^{-1}$ are periodic, $\|\Phi\|$ and $\|\Phi^{-1}\|$ are both bounded. Then, it can be shown that there exist $\bar \gamma_{k\ell}>0,k,\ell=1,\dots,r$, such that 
\begin{align*}
	\|\Phi^{-1} N\PC{k}(y) \Phi\|\le \sum_{\ell=1}^{r}\bar \gamma_{k\ell} \|z_\ell\|
\end{align*}
for each $k$ (see Lemma~\ref{lemma:bound:pert} in Appendix~\ref{proof:general} for more details). 

\begin{theorem}[{Sufficient condition for vibrational stabilization}]\label{theorem:general}
	Assume that $\bar J=\blkdiag(\bar J\PC{1},\dots,\bar J\PC{r})$ in Eq.~\eqref{J-bar} is Hurwitz. Let $\bar X_k$ be the solution to the Lyapunov equation 
	\begin{align}\label{Ly:controlled}
		{(\bar J\PC{k})^\top \bar X_k+\bar X_k \bar J\PC{k}=-I}.
	\end{align}
	Define the matrix $S=[s_{k \ell}]_{r\times r} $ with
	\begin{align}\label{matrix:S}
		s_{k\ell}=\Big\{\begin{matrix*}[l]
			\lambda_{\max}^{-1}(\bar X_k)-\bar\gamma_{k k} , &\text{ if }k=\ell,\\
			-\bar \gamma_{k \ell}, &\text{ if }k \neq\ell.
		\end{matrix*}
	\end{align}
	If $S$ is an $M$-matrix, then  there exists $\varepsilon^*>0$ such that, for any $\varepsilon<\varepsilon^*$:
	
	(i) the equilibrium $x=0$ of the system \eqref{compact:main} is exponentially stable uniformly in $y$; 
	
	(ii) the cluster synchronization manifold $\CM$ of the system \eqref{Kuramoto:controlled} is exponentially stable. 
\end{theorem}

The proof can be found in Appendix~\ref{proof:general}. Note that a similar theorem was presented in our conference paper \cite{YQ-DSB-FP:22a}. Applying a complete instead of a partial averaging technique in \eqref{change:time:scale}, we obtain a tighter condition by identifying smaller $\bar \gamma_{k\ell}$ in \eqref{matrix:S}. 
Theorem~\ref{general} provides a sufficient condition for vibrational control inputs to stabilize the cluster synchronization. To design an effective vibrational control law that stabilizes $\CM$, one just need to ensure that vibrations satisfy the following three conditions: 
i) $\bar J$ in \eqref{J-bar} is Hurwitz, 
ii) $S$ defined in \eqref{matrix:S} is an $M$-matrix, and
iii) the frequency of the vibrations is sufficiently high, i.e., $\varepsilon>0$ is sufficiently small.


\noindent
\textbf{Connections with robustness of linear systems:} Consider a stable linear system 
\begin{equation*}
	\dot x =Ax.
\end{equation*}
Some earlier works (e.g., \cite{PRV-TM:80,YR:85}) use 
\begin{equation}\label{measure:robust}
	\CR(A):=\lambda^{-1}_{\max}(X)
\end{equation}
to measure its robustness,  where $X$ is the solution to the Lyapunov equation
\begin{equation*}
	A^\top X +XA =-I. 
\end{equation*}
A larger $\CR(A)$ means that the system is more robust.

In our case, from \eqref{linearized}, the uncontrolled intra-dynamics around the manifold $\CM$ are described by 
\begin{align}\label{uncontrolled}
	&\dot x_k=J\PC{k} x_k+f\PC{k}_\inter(x,y),&k=1,\dots,r,
\end{align}
where $J\PC{k}$ is stable and $f\PC{k}_\inter(x,y)$ is taken as the vanishing perturbation.Here, $x_k=0$ means synchronization of the oscillators in the $k$th cluster. Similarly, one can interpret that $\CR(J\PC{k})$ measures the robustness of synchronization in the $k$th cluster.  If the intra-cluster synchronization is sufficiently robust (i.e., $\CR(J\PC{k})$'s are large) to dominate the perturbations resulted from inter-cluster connections, the cluster synchronization is stable. A sufficient condition is constructed in \cite[Th. 3.2]{TM-GB-DSB-FP:18}. By contrast, if $\CR(J\PC{k})$'s are not large enough, the cluster synchronization can lose its stability. Yet, the robustness of the intra-cluster synchronization can be reshaped by introducing vibrations to the local network connections. The new robustness is instead measured by $\CR(\bar J\PC{k})$.   

Now, the question naturally arises: how to design vibrational control such that the robustness in the cluster can be improved? We aim to provide answers in the next section.

\section{Improving Robustness by Vibrational Control}\label{sec:robustness}

The primary objective of this section is to demonstrate the design of vibrational control with the intention of enhancing the robustness of synchronization within each cluster. As observed in the concluding part of the previous section, the robustness is intimately linked to the linearized system. Hence, we commence by examining the linear system and subsequently explore the applicability of the findings from linear systems to Kuramoto-oscillator networks.

\subsection{Linear Systems}\label{subsec:linear}
Given a linear system 
\begin{equation}\label{linear}
	\dot x =A x,
\end{equation}
where $x\in\R^n, A\in\R^{n\times n}$, and $A$ is assumed to be Hurwitz. Consider a control matrix $V(t)$ that influences the system parameters in $A$, resulting in the following controlled system
\begin{equation}\label{controlled_net_compact}
	\dot x = \big(A+\frac{1}{\varepsilon}V\big(\frac{t}{\varepsilon}\big) \big) x,
\end{equation} 
where $\varepsilon>0$ is a small constant that determines the frequencies of the vibrations.


Vibrational control can improve the robustness of a stable system. To show this, we follow similar step as in Section~\ref{results:general} to associate \eqref{controlled_net_compact} with the averaged system 
\begin{equation}\label{linear:averaged}
	\dot{\bar x} = \bar A \bar x, 
\end{equation}
where
\begin{equation*}
	\bar A =\lim_{T\to \infty}\frac{1}{T} \int_{t=0}^{T} \Phi^{-1}(t,t_0) A \Phi(t,t_0) d t,
\end{equation*}
and $\Phi(t,t_0)$ is the state transition matrix of the system
\begin{equation*}
	\dot z = V(t) z. 
\end{equation*}

When $\bar A$ is Hurwitz, the controlled system \eqref{controlled_net_compact} behaves like \eqref{linear:averaged} on average. Then, one can interpret that vibrational control changes the system matrix from $A$ to $\bar A$ in an in-average sense. Vibrational inputs can be design to carefully modify the elements in $A$ so that $\CR(\bar A)$ is larger than $\CR(A)$, improving robustness. 

However, which elements in $A$ and how they can be changed is a challenging problem. An earlier attempt has been made in \cite{SMM:80}. Here, we aim to generalize their result by utilizing some graph-theoretical approaches. 

Specifically, we associate the uncontrolled system $\dot x =Ax$ with a weighted directed network $\CG_A=(\CV,\CE_A,A)$. Here, $\CV=\{1,2,\dots\}$, and there is a directed edge from $i$ to $j$, i.e., $(i,j)\in\CE_A$ if $a_{ji}\neq 0$. The matrix $A$ just becomes the weighted adjacency matrix. Likewise, one can associate the averaged controlled system \eqref{linear:averaged} with a weighted directed network $\CG_{\bar A}=(\CV,\CE_{\bar A},{\bar A})$, which we refer to as the \textit{functioning network}. 
Then, changing elements in $A$ reduces to alter the weights in the network $\CG_A$ by vibrational control. 

\begin{definition}
	The edge $(i,j)\in \CE_A$ is said to be vibrationally increasable if there exists a vibrational control $V(t)$ such that the weight of $(i,j)$ is increased $\CE_{\bar A}$, i.e.,  $\bar a_{ji}>a_{ji}$. It is said to be vibrationally decreasable if there exists a vibrational control $V(t)$ such that the weight of $(i,j)$ is decreased $\CE_{\bar A}$, i.e.,  $\bar a_{ji}<a_{ji}$.
\end{definition}

\begin{lemma}\label{lemma:inc_dec}
	Consider an edge $(i,j)\in \CE_A$, where $i\neq j$. It is vibrationally increasable if there is an edge in the reverse direction that has a negative weight, i.e., $a_{ij}<0$. It is vibrationally decreasable if there is an edge in the reverse direction that has a positive weight, i.e., $a_{ij}>0$. 
\end{lemma}

 \begin{figure}[t]
	\centering
	\includegraphics[scale=1]{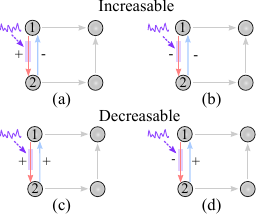}
	\caption{Illustration of vibrationally increasable and decreasable edges. (a) and (b): The red edge (i.e., $(1,2)$) is vibraionally increasable since the edge in the reverse direction has negative weight. (c) and (d): It is vibrational decreasable since the reverse edge is positively-weighted. Injecting a vibration to the red edge itself can functionally increase or decrease its weight in the corresponding cases.}
	\label{illus:decr_incr}
\end{figure}

\begin{figure*}[t]
	\centering
	\includegraphics[scale=1]{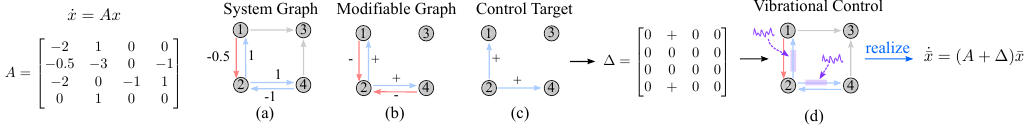}
	\caption{Illustration of the method to improve robustness of a linear stable system. (a) The directed graph associated with $\dot x =A x$. (b) Modifiable graph, where edges can be vibrationally changed (the signs indicate whether they can be increased or decreased). (c)-(d) If a matrix $\Delta$ corresponds to a singed and directed graph that is  directed acyclic and (ii) is a subgraph of the modifiable graph, then there exists vibrational control to realize the averaged system $\dot {\bar x} =(A+\Delta)\bar x$. Careful design of $\Delta$ can increase the robustness of the original system.}
	\label{theo:linear}
\end{figure*} 

The proof of this lemma can be found in Appendix~\ref{pf:inc_dec}.
An illustration of vibrationally increasable and decreasable edges can be found in Fig.~\ref{illus:decr_incr}. When an edge satisfies the corresponding conditions, we further find that directly injecting a vibration to it can functionally increase or decrease its weight.
Note that the conditions we identified here are just sufficient ones. Edges that do not satisfy these conditions may be also vibrationally increasable and decreasable.

Yet, we restrict our attention to the edges that satisfy the conditions in Lemma~\ref{lemma:inc_dec}. Then, based on them, we define two sets $\CE_{\rm inc}:=\{(i,j)\in \CE_A: a_{ij}<0\}$, and $\CE_{\rm dec}:=\{(i,j)\in \CE_A: a_{ij}>0\}$, which are vibrationally increasable and decreasable edges, respectively. Subsequently, we define a directed and signed graph $\CG^{\rm mod}_A:=(\CV,\CE_{\rm mod},S)$, where $\CE_{\rm mod}=\CE_{\rm inc} \cup  \CE_{\rm dec}$, and $S=\diag([\sign(a_{ij})]_{(i,j)\in \CE_{\rm mod}})$. We refer to $\CG^{\rm mod}_A$ as the \textit{modifiable graph} of $\CG_A$. An example is shown in Fig.~\ref{theo:linear} (b).

\begin{theorem}\label{vibrational:design}
	Consider a matrix $\Delta=[d_{ij}]\in\R^{n\times n}$, and let $\CG_\Delta:=(\CV,\CE_\Delta,S_\Delta)$ be the directed and signed graph associated with it. If the following conditions are satisfied:
	
	(i) $\CG_\Delta$ is a directed acyclic graph and $\CG_\Delta \subset \CG_A^{\rm mod}$,
	
	(ii) $A+\Delta$ is Hurwitz,
	
	\noindent
	Then, the system matrix $\bar A$ of \eqref{linear:averaged} becomes $\bar A=A+\Delta$ if the vibrational control inputs injected to the edges in $\CG_\Delta$ are
	\begin{equation}\label{imp_robu:linear}
		v_{ij}(t)=\frac{1}{\varepsilon}\sqrt{\frac{-2d_{ij}}{a_{ji}}} \sin \Big(\frac{\beta_{ij}t}{\varepsilon}\Big),
	\end{equation}
	where $\beta_{ij}$'s are incommensurable\footnote{Two non-zero real numbers $a$ and $b$ are said to be incommensurable if their ratio ${a}/{b}$ is not a  rational number}.
\end{theorem}

This theorem provides a method to improve the robustness of a system described by \eqref{linear}. As illustrated in Fig.~\ref{theo:linear}, if one can choose a matrix $\Delta$ satisfying the conditions (i) and (ii) and $\CR(A+\Delta)>\CR(A)$, the vibrational control inputs in \eqref{imp_robu:linear} improves the robustness of the system from $\CR(A)$ to $\CR(A+\Delta)$.  We wish to mention that it is likely impossible to improve the system to any desired robustness by just adding a matrix $\Delta$, especially when $\Delta$ is constrained by the graph structure. There are some interesting open questions. For instance, what is the realizable range of robustness levels? How to design $\Delta$ and the subsequent vibrational control to realize a desired and reasonable robustness? 

\subsection{Kuramoto-Oscillator Networks}\label{subsec:kuramo}
Observe that, in each cluster, oscillators have an identical frequency, and each pair is coupled by bidirected edges with asymmetric strengths. To study how vibrational control can improve the robustness of synchronization in each cluster, we consider the following homogeneous Kuramoto model:
\begin{equation}
	\dot \phi_i = \beta + \sum_{j=1}^{n} w_{ij} \sin(\phi_j -\phi_i),
\end{equation}
where $i=1,\dots,n$, and $w_{ij}$'s describes the directed network $\BG=(\BV,\BE)$ with $|\BE|=m$. Let $B\in \R^{n \times m}$ be the incidence matrix of $\BG$.  Select a directed spanning tree $\BG_{\rm span}$ in $\BG$, and let $\hat B\in \R^{n\times (n-1)}$ be its incidence matrix. Denote $\phi=[\phi_1,\dots,\phi_m]$ and $W=\diag([w_{ij}]_{(i,j)\in \BE})\in\R^{m\times m}$.

Let $x= \hat B^\top \phi\in \R^{n-1}$, and following similar steps as Appendix~\ref{derivatin} one can derive that
\begin{equation}\label{homog:compact}
	\dot x = - \hat B^\top \BB W \sin (R_1 x),
\end{equation}
where $\BB$ is obtained by replacing the positive elements in the incidence matrix $B$ by $0$ and
\begin{equation*}
	R_1 = \begin{bmatrix}
		B^\top \hat B^\dagger\\
		-B^\top \hat B^\dagger
	\end{bmatrix}.
\end{equation*}

With vibrations injected into the edges in the network, we have the control model
\begin{equation}\label{homog:compact:controlled}
	\dot x = - \hat B^\top \BB \Big(W+\frac{1}{\varepsilon}V\big(\frac{t}{\varepsilon}\big)\Big) \sin (R_1 x),
\end{equation}
where $V(t)=\diag([v_{ij}]_{(j,i)\in\BE})\R^{m\times m}$. 
Linearizing the system \eqref{homog:compact:controlled} at $x=0$, we obtain
\begin{equation}\label{ctr_KM:linearized}
	\dot x = - \hat B^\top \BB \Big(W+\frac{1}{\varepsilon}V\big(\frac{t}{\varepsilon}\big)\Big) R_1 x. 
\end{equation}
Denote $J=-\hat B^\top \BB W R_1$. Similar to the previous section, one can associate \eqref{ctr_KM:linearized} with following averaged system
\begin{equation}
	\dot{\bar x} = \bar J \bar x
\end{equation}
where 
\begin{equation*}
	\bar J =\lim_{T\to \infty}\frac{1}{T} \int_{t=0}^{T} \Phi^{-1}(t,t_0)J\Phi(t,t_0) d t,
\end{equation*}
and $\Phi(t,t_0)$ is the state transition matrix of the system
\begin{equation}
	\dot z = P(t) z, \text{where } P(t) =  - \hat B^\top \BB V(t)R_1.
\end{equation}

Due to the presence of the matrices $\hat B,\BB$, and $R_1$, $P(t)$ has very complex dependence on the vibrations injected to the edges in the Kuramoto-oscillator network. As a consequence, it becomes very challenging to design vibrational control $V(t)$ to modify the elements in $J$ so that the robustness of synchronization can be improved.

\begin{figure*}[t]
	\centering
	\includegraphics[scale=1]{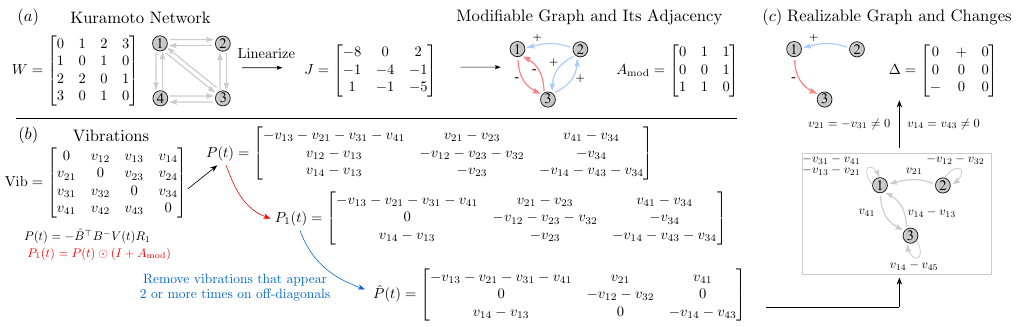}
	\caption{Illustration of making predictable changes in the Kuramoto-oscillator network. (a) One can linearize the Kuramoto model and obtain a linear system $\dot x = J x$. Following the steps in Section~\ref{sec:robustness}-\ref{subsec:linear}, a modifiable graph and its adjacency matrix can be found. (b) We investigate how vibrations in the Kuramoto network influence the linearized system by computing $P=\hat B^\top \BB V(t)R_1$. From the modifiable graph in (a), we only keep elements in $P(t)$ that corresponds to modifiable edges, resulting in $P_1(t)$. Further, to make the design of vibrational tractable, we remove the vibrations that appear two or more times in off-diagonal positions and obtain $\hat P(t)$. (c) We associate $\hat P(t)$ with a directed graph $\CG_{\hat P}$ that allows for self-loops. Then, it remains to configure vibrations such that the directed graph has no directed cycles (including self-loops). For instance, one can set $v_{21}(t)=-v_{31}(t)$, $v_{14}(t)=v_{45}(t)$, and any other vibrations to zero to realize the changes indicated in the upper panel. Any change that has the same pattern as $\Delta$ can by realized by carefully choosing the amplitudes of vibrations $v_{21}(t)$ and $v_{14}$.}
	\label{design:kuramoto}
\end{figure*}

The following lemmas provides a necessary condition such that the robustness can be improved.

\begin{lemma}[\textit{Necessary condition for robustness improvement}]\label{necessary}
	There exists a vibrational control input such that $\CR(\bar J)> \CR(J)$ \emph{only if} the following two conditions are both satisfied
	
	(i) $\BG$ is a bidirected complete graph\footnote{A bidirected complete graph is a graph where each node receives a connection from any other node.};
	
	(ii) for each $i$, $w_{ij}$'s are identical for all $k\neq j$. 
\end{lemma} 

In this section, we are interested in the situation where the robustness of the system can be shaped. Our goal is to provide a tractable and predictable approach to design vibrational control to functionally modify the elements in $J$, aiming to improving its robustness. The results we established for linear systems in Section~\ref{sec:robustness}-\ref{subsec:linear} will be used. 

Specifically, the methods consists of the following steps (an example if provided in Fig.~\ref{design:kuramoto} to illustrate the procedure). 

(1) First, we compute the Jacobian matrix $J$ and associate it with a weighted directed graph $\BG_J$. 

(2) Following the same steps as those in Section~\ref{sec:robustness}-\ref{subsec:linear}, one can identify a modifiable graph from $\BG_J$, indicating which edges in $\CG_J$ are increasable and decreasable. Let $A_{\rm mod}$ be the unweighted adjacency matrix of $\BG_J$.

Different from the linear systems in the previous section, the edges can not be directly modified since each change has to be induced by vibrating the connections in the original Kuramoto-oscillator network. 

(3) Compute $P(t)=-\hat B^\top \BB V(t)R_1$, which captures how vibrations injected to the edges in the Kuramoto-oscillator network affects the edges in $\BG_J$. Since whether an edge in $\BG_J$ can be altered is determined by its modifiable graph, we let $P_1(t) = P(t)\odot (I+A_{\rm mod})$, which only keeps the elements that corresponds to modifiable edges. 

 While configuring control inputs, one needs to deal with the situation that a vibration introduced to a single edge in the Kuramoto-oscillator network can bring changes to multiple edges in $\BG_J$. Therefore, one often needs to combine multiple vibrations to avoid that.

(4) To make the design more analytically tractable, we remove the vibrations that appear two or more times in the off-diagonal positions in $P_1(t)$, and obtain $\hat P(t)$. 

(5) We can associate $\hat P(t)$ with a directed graph $\BG_{\hat P}$. Now, one can configure vibrational control inputs such that $\BG_{\hat P}$ does not contain a directed cycle (including self-loops). Consequently, the resulting graph determines realizable changes to $J$ that vibrations can bring in, which we refer to as a \textit{realizable graph}. For any $\Delta$, there always exists a vibrational control that functionally changes $J$ to $\bar J =J+\Delta$ if the associated directed and sign graph of $\Delta$ is the same as a realizable graph (see Fig.~\ref{design:kuramoto}~(c)). One can simply use the results in Theorem~\ref{vibrational:design} to design vibrational inputs. 

\begin{remark}
	We remark that the vibrations that only appear in the diagonal positions of $P(t)$ play an important role in the procedure of design. They are often used to cancel other vibrations' influence on the diagonal positions, ensuring they can only affect one off-diagonal element $P(t)$. To ensure as many vibrations as possible to appear only  in the diagonal positions, one way is to select a spanning tree as short as possible to define $\hat B$. For instance, in Fig.~\ref{design:kuramoto} (a), one can choose the spanning tree consisting of $\{(3,1),(3,2),(3,4)\}$. 
	\end{remark}

\subsection{Design of Vibrational Control for Cluster Synchronization Stabilization}
Now, we can use the results in the previous section to design vibrations to stabilize cluster synchronization. 

Recall that oscillators in each cluster have an identical frequency. Following the same procedure as in the previous section, one can identify a realizable graph for each cluster, defining changes to each Jacobian matrix that can be realized by vibrational control. Denote the realizable graphs by $\BG^{(1)}, \BG^{(2)}, \dots,\BG^{(r)}$, and the Jacobian matrices by $J^{(1)},J^{(2)},\dots,J^{(r)}$. Let $\Delta^{(k)}\in \R^{n_k \times n_k}$ be a matrix constrained by $\BG\PC{k}$ for $k=1,2,\dots,r$. 

\begin{corollary}\label{coro}
	Assume that there exist matrices $\Delta^{(k)}\in \R^{n_k \times n_k}, k=1,2,\dots,r$ that satisfy the following conditions:
	
	(i) Each $\Delta^{(k)}\in \R^{n_k \times n_k}$ is constrained by a realizable graph $\BG\PC{k}$ for $k=1,2,\dots,r$. 
	
	(ii) The matrix $S=[s_{k \ell}]_{r\times r} $ defined by
	\begin{align*}
		s_{k\ell}=\Big\{\begin{matrix*}[l]
			\CR(J\PC{k}+\Delta\PC{k})-\bar\gamma_{k k} , &\text{ if }k=\ell,\\
			-\bar \gamma_{k \ell}, &\text{ if }k \neq\ell.
		\end{matrix*}
	\end{align*}
	is an $M$-matrix.
	
	Then, the vibrational control inputs that functionally change $J\PC{k}$ to $J\PC{k}+\Delta\PC{k}$ stabilize the cluster synchronization manifold $\CM$. 
\end{corollary}

We wish to mention that one can follow the same steps as in Sections~\ref{sec:robustness}-\ref{subsec:linear} and \ref{subsec:kuramo} to design the amplitudes and frequencies of vibrational inputs. 

\section{Numerical Study}\label{numerical}
In this section, we employ an example to show how to design vibrations to stabilize cluster synchronization in a Kuramoto-oscillator network. 

\begin{figure}[t]
	\centering
	\includegraphics[scale=1]{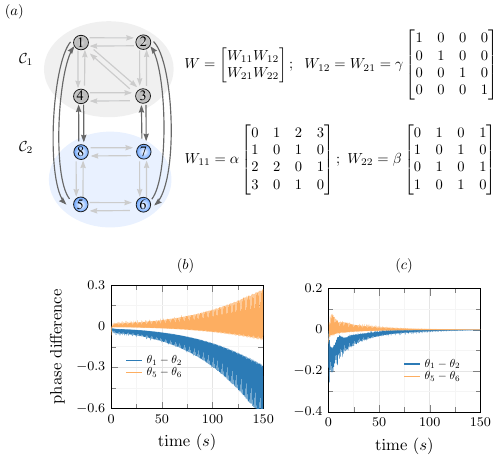}
	\caption{Vibrational stabilization of a cluster synchronization manifold. (a) The network structure and the connection weights, where $\alpha=0.05$, $\beta=1$, and $\gamma=3$ (b) Phase differences without control, indicating that the cluster synchronization is unstable. (c) Phase differences under vibrational control to the cluster $\mathcal{C}_1$. showing that the cluster synchronization has been stabilized by local vibrations. The natural frequencies in $\mathcal{C}_1$ and $\mathcal{C}_2$ are $\omega_1=1$ and $\omega_2=10$, respectively.}
	\label{fig:kuramoto_stabilization}
\end{figure}

The network we consider is shown in Fig.~\ref{fig:kuramoto_stabilization}~(a). Partitioning the network into two clusters $\CC_1$ and $\CC_2$, Assumption~\ref{invariance} is satisfied so that the corresponding cluster synchronization manifold $\CM$ is invariant. However, this pattern of cluster synchronization is unstable (see in Fig.~\ref{fig:kuramoto_stabilization}~(b)). Then, we want to design a vibrational control to stabilize it. 

We observe that the first cluster has the same network structure as that in Fig.~\ref{design:kuramoto}. Within each cluster, one can derive that the linearized system $\dot x\PC{k}= J^{(k)} x\PC{k}$ has
\begin{equation*}
	J^{(1)}=\alpha\begin{bmatrix}
		-8 & 0 & 2\\
		-1 & -4 & -1\\
		1 & -1 & -5\\
	\end{bmatrix}, \quad
	J^{(2)}=\begin{bmatrix}
		-3 & 0 & 1\\
		-1 & -2 & 1\\
		1 & 0 & -3\\
	\end{bmatrix}.
\end{equation*}

By some simple computation, one notices that the robustness of synchronization within $\CC_1$ is small. Fig.~\ref{design:kuramoto} has identified a realizable graph and changes, and we use it to modify the elements in $J^{(1)}$. Particularly, we choose 
\begin{equation*}
	\Delta\PC{1}=0.05 \begin{bmatrix}
		0 & 1& 0\\
		0 &0 &0\\
		-1 &0 & 0
	\end{bmatrix}.
\end{equation*}
One  can compute this change improves the robustness from $\CR(J^{(1)})=0.305$ to $\CR(J^{(1)}+\Delta\PC{1})=0.332$.  

 Following Theorem~\ref{vibrational:design}, we inject the vibrations below to realize these changes:

\begin{align*}
	&a_{21}(t)=-a_{31} (t) = \frac{k_1}{\varepsilon}\sin\big(\frac{\beta_1}{\varepsilon}\big),\\
	&a_{14}(t)=a_{45} (t) = \frac{k_2}{\varepsilon}\sin\big(\frac{\beta_2}{\varepsilon}\big),
\end{align*}
where $\beta_1=1$, $\beta_2=\sqrt{2}$, and
\begin{equation*}
	k_1=\sqrt{-\frac{0.05}{J\PC{1}_{21}}}, \text{and } k_2= \sqrt{\frac{0.05}{J\PC{1}_{13}}}.
\end{equation*}

Let $\varepsilon=0.01$, and one can observe from Fig.~\ref{fig:kuramoto_stabilization} that the cluster synchronization is stabilized. We wish to mention that the condition in Theorem~\ref{theorem:general} and Corollary~\ref{coro} are not even satisfied. This indicates that the condition we have identified is still a bit conservative.  More tight conditions call for future studies. However, it is worth emphasizing the power of vibrational control since a slight improvement on the robustness effectively stabilizes the cluster synchronization.

\appendix

\subsection{Derivation of the Compact System}\label{derivatin}

Let 
\begin{equation}
	W=\begin{bmatrix}
			W_\intra&0\\
			0&W_\inter
		\end{bmatrix}
\end{equation}
where $W_\intra:=\diag\{a_{ij},(i,j)\in \CG_\intra\}$ and $W_\inter:=\diag\{a_{ij},(i,j)\in \CG_\inter\}$ are diagonal weight matrices of $ \CG_\intra$ and $ \CG_\inter$, respectively. Similarly, one can use 
\begin{equation}
	V(t)=\begin{bmatrix}
		V_\intra(t)&0\\
		0&V_\inter(t)
	\end{bmatrix}
\end{equation}
to denote the vibrations injected to the edges that corresponds to $W$. Then, one can rewrite the controlled system \eqref{Kuramoto:controlled} into
\begin{equation}
	\dot \theta = \omega - \BB (W+V(t)) \sin(B^\top \theta), 
\end{equation}
where $\BB$ is obtained by replacing the positive elements in the incidence matrix $B$ by $0$.
Recall that $x=\hat B_\intra^\top \theta$ and $y=\hat B_\inter^\top \theta$. Then, it holds that 
\begin{subequations}\label{derive:intermediate}
\begin{align}
	&\dot x = -  \hat B_\intra^\top \BB (W+V(t)) \sin(B^\top \theta),\\
	&\dot y = \hat B_\inter^\top \omega - \hat B_\inter^\top \BB (W+V(t)) \sin(B^\top \theta), 
\end{align}
\end{subequations}
where the fact that intra-cluster natural frequency difference are zero has been used.
Now, it remains to write $B^\top \theta$ into a function of $x$ and $y$. 

To this end, we provide the following instrumental lemma. Below, $\bar B$ is the incidence matrix of the unidirected counterpart of $\CG$ (removing one edge from each bidirected edge in $\CG$); and $\bar B_\inter$ and $\bar B_\intra$ are also the unidirected counterpart of $ B_\inter$ and $ B_\intra$, respectively. 
\begin{lemma}\label{incidence:transfer}
	For the incidence matrices $B$ and $\hat B$, there exists 
	\begin{equation*}
		R=\scalemath{0.85}{\begin{bmatrix}
			R'_1& 0 \\
			-R'_1 & 0\\
			R'_2 & R'_3\\
			-R'_2 & -R'_3
		\end{bmatrix}}:= \begin{bmatrix}
		R_1 & 0\\
		R_2 &R_3
		\end{bmatrix}
	\end{equation*}
	such that $B^\top =R\hat B^\top$,  where 
	\begin{align*}
		&R'_1=\scalemath{0.85}{\begin{bmatrix}
			(\bar B^{(1)}_{\intra})^\top ((\hat B^{(1)}_{\intra})^\top)^\dagger & \cdots &0\\
			\vdots & \ddots & \vdots\\
			0& \cdots & (\bar B^{(1)}_{\intra})^\top ((\hat B^{(1)}_{\intra})^\top)^\dagger
		\end{bmatrix}} \\
		&R'_2=\bar B_\inter^\top (\hat B_\intra ^\top P_\inter)^\dagger, \text{and } R'_3= \bar B_\inter^\top (\hat B_\inter P_\intra)^\dagger
	\end{align*}
	with $\SmallMatrixNBla{P_\intra=I_n-\hat B_\intra \hat B_\intra^\dagger}$ and $\SmallMatrixNBla{P_\inter=I_n-\hat B_\inter \hat B_\inter^\dagger}$. \QEDA
\end{lemma}
The proof of this lemma can be straightforwardly adapted  from Lemma~5 in \cite{YQ-DSB-FP:22a} and from \cite{KR-IH:2021}, which is omitted here.  

Applying Lemma~\ref{incidence:transfer} to Eqs. \eqref{derive:intermediate}, one can derive the compact system \eqref{compact_form}, where
\begin{subequations}\label{expr:function}
\begin{align}
	&f_\intra (x) = -\hat B_\intra^\top \BB_\intra W_\intra \sin (R_1x) \\
	&f_\inter (x,y)= - \hat B_\intra^\top \BB_\inter W_\inter \sin(R_2 x +R_3y) \label{expr:inter-pert}\\
	&f_{\ctl}(x,y,V(t))=- \hat B_\intra^\top \BB_\intra V_\intra(t) \sin (R_1x) \\
	&- \hat B_\intra^\top \BB_\inter V_\inter(t) \sin(R_2 x +R_3y), \\
	&g(x,y)=\hat B_\inter^\top \omega - \hat B_\inter^\top \BB_\intra W_\intra \sin (R_1x) \\
	&- \hat B_\inter^\top \BB_\inter W_\inter \sin(R_2 x +R_3y) \\
	&g_\ctl(x,y,V(t))=- \hat B_\inter^\top \BB_\intra V_\intra(t) \sin (R_1x) \\
	&- \hat B_\inter^\top \BB_\inter V_\inter(t) \sin(R_2 x +R_3y),
\end{align}
\end{subequations}

\subsection{Proof of Lemma~\ref{lemma:linearized}}\label{proof:linearized}
Since $x=0$ is exponentially stable uniformly in $y$ for the system \eqref{linearized}, according to the converse Lyapunov theorem  (see \cite[Th. 4.4]{haddad2011nonlinear} and \cite{QY-KY-BDOA-CM:2021}) there exists $\CD=\{x\in\R^{n-r}:\|x\|\le \rho_1\}$ and a  continuously
differentiable function function $V:[0,\infty]\times \CD\times \R^r\to \R $ such that 
\begin{align*}
			\frac{\partial V}{\partial t}+\frac{\partial V}{\partial x}\big( (J + P(t))x+N(x,y)x\big)\le -c_1\|x\|^2
\end{align*}
and $\|\frac{\partial V}{\partial x}\|\le c_2\|x\|$ for some constants $c_1,c_2> 0$.  
Let 
\begin{equation*}
	h(t,x)=f_\intra(x)+f_\ctr(U(t),x)+f_\inter(x,y)
\end{equation*}

and 
\begin{equation*}
\Delta(t,x)=h(t,x)-(J + P(t))x-N(x,y).
\end{equation*}
It can be checked that $\partial h/\partial x$ is bounded and Lipschitz on $\CD$. Then, similar to the proof of \cite[Th. 4.13]{HKK:02-bis}, one can show that $\|\Delta(t,x)\|\le c_3 \|x\|^2$ for some $c_3>0$. The time derivative along the system  \eqref{compact_form} satisfies
\begin{align*}
	&{\frac{\partial V}{\partial t}+\frac{\partial V}{\partial x}\big( (J + P(t))x+\Delta(t,x)\big)}\\
	&\le -c_1\|x\|^2+c_2c_3\|x\|^3\\
	&\le -(c_1-c_2c_3\rho)\|x\|^2, \forall \|x\|< \rho.
\end{align*}
Choosing $\rho=\min\{\rho_1,c_1/c_2c_3\}$ completes the proof.

\subsection{Proof of Theorem~\ref{theorem:general}} \label{proof:general}
One can observe that (i) implies (ii). Then, it suffices to prove the exponential stability of $x=0$ for \eqref{change:time:scale}. To do that, we first present the following lemma, whose proof follows similar lines as Lemma 3.1 of \cite{TM-GB-DSB-FP:18}.
\begin{lemma}[Growth bound of perturbations]\label{lemma:bound:pert}
	There exist some constants $\bar \gamma_{k\ell}>0$, $k,\ell=1,\dots,r$, such that, for any $k$, it holds that 
	$
	\|\Phi^{-1} N\PC{k}(y) \Phi\|\le \sum\nolimits_{\ell=1}^{r}\bar \gamma_{k\ell} \|z_\ell\|.
	$\QEDA
\end{lemma}

Let $\scalemath{0.85}{V_k=z_k^\top \bar X_kz_k}$, and we have 
\begin{equation*}
	{\|\partial V_k/\partial z_k \| \le \lambda_{\max}( \bar X_k)\|z_k\|}.
\end{equation*}
	 Choose $\scalemath{0.85}{V(z)= \sum\nolimits_{k=1}^{r}d_k V_k}$ as a Lyapunov candidate. The time derivative of $V(z)$ satisfies
\begin{align*}
	{\dot V(z)}&=\sum_{k=1}^{r} d_k[z_k^\top((\bar J\PC{k})^\top  \bar X_k+  \bar X_k \bar J\PC{k})z_k \\
		&+\frac{\partial V}{\partial z_k} \Phi^{-1} N(y) \Phi]\\
	&{\le \sum_{k=1}^{r} d_k[-\|z_k\|^2 + \lambda_{\max}( \bar X_k)\sum_{\ell =1}^{r}\bar \gamma_{k \ell} \|z_k\|\|z_\ell\| ]},
\end{align*}
where the second inequality has used Lemma~\ref{lemma:bound:pert}.

Let ${D:=\diag(d_1,\dots,d_r)}$ and ${\hat S=[\hat s_{ij}]_{r\times r}}$ where 
\begin{align*}
	\hat s_{k\ell}=\Big\{\begin{matrix*}[l]
		1- \lambda_{\max}( \bar X_k) \bar \gamma_{k k} , &\text{ if }k=\ell,\\
		-\lambda_{\max}( \bar X_k)\bar \gamma_{k \ell}, &\text{ if }k \neq\ell.
	\end{matrix*}
\end{align*}
Then, one can rewrite ${\dot V(z)\le -\frac{1}{2}z^\top (DS+S^\top D)z}$. By assumption, $S$ is an $M$-matrix, and so is $\hat S$ since $\hat S=S \cdot\diag(\lambda_{\max}( \bar X_1),\dots,\lambda_{\max}( \bar X_r))$. It follows from \cite[Th. 9.2]{HKK:02-bis} that the system \eqref{average} is exponentially stable.  

Following similar steps as in \cite[Th. 10.4]{HKK:02-bis}, one can prove that there exists $\varepsilon^*>0$ such that for any $\varepsilon<\varepsilon^*$, $z=0$ is exponentially stable uniformly in $y$ for the system \eqref{coordinated}. Since $x(s)=\Phi(s,s_0)z(s)$ and $\|\Phi\|$ is bounded, then $x=0$ is also exponentially stable uniformly in $y$ for \eqref{change:time:scale}, which completes the proof. 

\subsection{Proof of Lemma~\ref{lemma:inc_dec}} \label{pf:inc_dec}

	To construct the proof, one just needs to find vibrational control inputs that increase/decrease the weight of the edge $(i,j)$ functionally in both situations.
	
	Consider a vibrational control that is only injected to the edge $(i,j)$, i.e., $V(t)$ in \eqref{controlled_net_compact} satisfies $v_{pq}(t)=0$ for any $(p,q)\neq (i,j)$ and $V_{ji}(t)\neq 0$. One can label the nodes in the expanded network $\CG$ such that $i=1$ and $j=2$. Then, the vibrational control matrix becomes 
	\begin{equation*}
		V(t)=\begin{bmatrix}
			0&0&\cdots&0 \\
			v_{21}(t)&0&\cdots&0\\
			\vdots&\vdots&\ddots&\vdots \\
			0&0&0&0 
		\end{bmatrix},
	\end{equation*}
	which has a quasi-lower-triangular form. Following the steps in \cite{SMM:80}, one can derive that $\bar A$ in the averaged system \eqref{linear:averaged} is
	\begin{equation*}
		\bar A = A + \bar B, \text{ where } \bar B = \begin{bmatrix}
			0&0&\cdots&0 \\
			b_{21}(t)&0&\cdots&0\\
			\vdots&\vdots&\ddots&\vdots \\
			0&0&0&0 
		\end{bmatrix},
	\end{equation*}
	where $b_{21} = -a_{12} \lim_{T \to \infty}\frac{1}{T} \int_{t=0}^{T} F^2_{21}(t) dt$ with $F_{21}(t)=\int_{0}^{t}v_{21}(\tau) d\tau$. If the edge from $1$ to $2$ has a positive weight, i.e., $a_{12}> 0$, $b_{21}<0$; $b_{21}>0$ if $a_{12}<0$, which completes the proof.
	
	\subsection{Proof of Theorem~\ref{vibrational:design}}
	If $\CG_\Delta$ is directed acyclic, according to \cite{JBJ-GG:00},  it can be topologically ordered. Therefore, one can  arrange the vertices of $\CG_\Delta$ as a linear ordering that is consistent with all edge directions. In other words, there exists a permutation matrix $P$ such that the matrix $\Delta'=:P\Delta P^{-1}$ is quasi-lower-triangular. One can let $A'=PAP^{-1}$. Let $x'=Px$, one can derive that
	\begin{align*}
		\dot x' = A'x'. 
	\end{align*}
	Now, to change $A$ to $A+\Delta$, it becomes to change $A'$ to $A'+ \Delta'$, with $\Delta'=[d'_{ij}]\in\R^{n\times n}$ being  quasi-lower-triangular. 
	
	Now, consider the following vibrations
	\begin{equation*}
		v_{ij}(t)=\frac{1}{\varepsilon}\sqrt{\frac{-2d'_{ij}}{a_{ji}}} \sin \Big(\frac{\beta_{ij}t}{\varepsilon}\Big),
	\end{equation*}
	for all pairs of $i$ and $j$ such that $d'_{ij}\neq 0$. 
	Following similar steps as in \cite{AMN-YQ-CAA-DSB-FP:22} and \cite{SMM:80}, one can derive that the system
	$
		\dot x'=(A'+V'(t))x'
	$
	has a stable averaged system $\dot x'=(A'+\Delta')x'$, which completes the proof.

\subsection{Proof of Lemma~\ref{necessary}}
When the conditions are satisfied, it holds that $J=-dI_{n-1}$, where $d=\sum_{j\neq i}w_{ij}+w_{ji}$. As a consequence, $\CR(\bar J)= \CR(J)$ for any vibrational control since
\begin{align*}
	\scalemath{0.9}{\bar J=\int_{0}^{2 \pi} \Phi^{-1}(t) J\Phi(t)=-d\int_{0}^{2 \pi} \Phi^{-1}(t) \Phi(t)dt}=J.
\end{align*}

\bibliographystyle{IEEEtran}
\bibliography{alias,FP,Main,New}

\begin{thebibliography}{10}
\providecommand{\url}[1]{#1}
\csname url@samestyle\endcsname
\providecommand{\newblock}{\relax}
\providecommand{\bibinfo}[2]{#2}
\providecommand{\BIBentrySTDinterwordspacing}{\spaceskip=0pt\relax}
\providecommand{\BIBentryALTinterwordstretchfactor}{4}
\providecommand{\BIBentryALTinterwordspacing}{\spaceskip=\fontdimen2\font plus
\BIBentryALTinterwordstretchfactor\fontdimen3\font minus
  \fontdimen4\font\relax}
\providecommand{\BIBforeignlanguage}[2]{{%
\expandafter\ifx\csname l@#1\endcsname\relax
\typeout{** WARNING: IEEEtran.bst: No hyphenation pattern has been}%
\typeout{** loaded for the language `#1'. Using the pattern for}%
\typeout{** the default language instead.}%
\else
\language=\csname l@#1\endcsname
\fi
#2}}
\providecommand{\BIBdecl}{\relax}
\BIBdecl

\bibitem{HG-PAA-DG-AA-KA:19}
G.~Hahn, A.~Ponce-Alvarez, G.~Deco, A.~Aertsen, and A.~Kumar, ``Portraits of
  communication in neuronal networks,'' \emph{Nature Reviews Neuroscience},
  vol.~20, no.~2, pp. 117--127, 2019.

\bibitem{FJ-AN:2011}
J.~Fell and N.~Axmacher, ``The role of phase synchronization in memory
  processes,'' \emph{Nature Reviews Neuroscience}, vol.~12, no.~2, pp.
  105--118, 2011.

\bibitem{Hammond2007}
C.~Hammond, H.~Bergman, and P.~Brown, ``Pathological synchronization in
  {P}arkinson's disease: Networks, models and treatments,'' \emph{Trends in
  Neurosciences}, vol.~30, no.~7, pp. 357--364, 2007.

\bibitem{JP-DCM-JJGR-SCA:2013}
P.~Jiruska, M.~De~Curtis \emph{et~al.}, ``Synchronization and desynchronization
  in epilepsy: Controversies and hypotheses,'' \emph{The Journal of
  Physiology}, vol. 591, no.~4, pp. 787--797, 2013.

\bibitem{REB-JB-SMM:86b}
R.~E. Bellman, J.~Bentsman, and S.~M. Meerkov, ``Vibrational control of
  nonlinear systems: Vibrational controllability and transient behavior,''
  \emph{IEEE Transactions on Automatic Control}, vol.~31, no.~8, pp. 717--724,
  1986.

\bibitem{BS-BTZ:97}
B.~Shapiro and B.~T. Zinn, ``High-frequency nonlinear vibrational control,''
  \emph{IEEE Transactions on Automatic Control}, vol.~42, no.~1, pp. 83--90,
  1997.

\bibitem{CX-TY-MI:2018}
X.~Cheng, Y.~Tan, and I.~Mareels, ``On robustness analysis of linear
  vibrational control systems,'' \emph{Automatica}, vol.~87, pp. 202--209,
  2018.

\bibitem{Krauss2021}
J.~K. Krauss, N.~Lipsman, T.~Aziz, A.~Boutet, P.~Brown, J.~W. Chang,
  B.~Davidson, W.~M. Grill, M.~I. Hariz, A.~Horn \emph{et~al.}, ``Technology of
  deep brain stimulation: current status and future directions,'' \emph{Nature
  Reviews Neurology}, vol.~17, no.~2, pp. 75--87, 2021.

\bibitem{TM-GB-DSB-FP:18}
T.~Menara, G.~Baggio, D.~S. Bassett, and F.~Pasqualetti, ``Stability conditions
  for cluster synchronization in networks of heterogeneous {K}uramoto
  oscillators,'' \emph{IEEE Transactions on Control of Network Systems},
  vol.~7, no.~1, pp. 302--314, 2020.

\bibitem{TM-GB-DSB-FP:19b}
------, ``A framework to control functional connectivity in the human brain,''
  in \emph{{IEEE} Conf.\ on Decision and Control}, Nice, France, Dec. 2019, pp.
  4697--4704.

\bibitem{TM-GB-DSB-FP:22}
------, ``Functional control of oscillator networks,'' \emph{Nature
  Communications}, vol.~13, p. 4721, 2022.

\bibitem{Pecora2014}
L.~M. Pecora, F.~Sorrentino, A.~M. Hagerstrom, T.~E. Murphy, and R.~Roy,
  ``Cluster synchronization and isolated desynchronization in complex networks
  with symmetries,'' \emph{Nature Communications}, vol.~5, no.~1, pp. 1--8,
  2014.

\bibitem{YSC-TN-AEM:17}
Y.~S. Cho, T.~Nishikawa, and A.~E. Motter, ``Stable chimeras and independently
  synchronizable clusters,'' \emph{Physical {R}eview {L}etters}, vol. 119,
  no.~8, p. 084101, 2017.

\bibitem{YQ-MC-BDOA-DSB-FP:20}
Y.~Qin, M.~Cao, B.~D.~O. Anderson, D.~S. Bassett, and F.~Pasqualetti,
  ``Mediated remote synchronization: the number of mediators matters,''
  \emph{IEEE Control Systems Letters}, vol.~5, no.~3, pp. 767--772, 2020.

\bibitem{EJ-SA-SJ-CJP-DRM:20}
J.~Emenheiser, A.~Salova, J.~Snyder, J.~P. Crutchfield, and R.~M. D'Souza,
  ``Network and phase symmetries reveal that amplitude dynamics stabilize
  decoupled oscillator clusters,'' \emph{arXiv preprint arXiv:2010.09131},
  2020.

\bibitem{Schaub2016}
M.~T. Schaub, N.~O'Clery, Y.~N. Billeh, J.-C. Delvenne, R.~Lambiotte, and
  M.~Barahona, ``Graph partitions and cluster synchronization in networks of
  oscillators,'' \emph{Chaos}, vol.~26, no.~9, p. 094821, 2016.

\bibitem{QY-KY-PO-CM:21}
Y.~Qin, Y.~Kawano, O.~Portoles, and M.~Cao, ``Partial phase cohesiveness in
  networks of networks of kuramoto oscillators,'' \emph{IEEE Transactions on
  Automatic Control}, vol.~66, no.~12, pp. 6100--6107, 2021.

\bibitem{FP-SA-MT:21}
P.~Feketa, A.~Schaum, and T.~Meurer, ``Stability of cluster formations in
  adaptive kuramoto networks,'' \emph{IFAC-PapersOnLine}, vol.~54, no.~9, pp.
  14--19, 2021.

\bibitem{FP-SA-MT:2020}
------, ``Synchronization and multicluster capabilities of oscillatory networks
  with adaptive coupling,'' \emph{IEEE Transactions on Automatic Control},
  vol.~66, no.~7, pp. 3084--3096, 2020.

\bibitem{KR-IH:2021}
R.~Kato and H.~Ishii, ``Averaging and cluster synchronization of {K}uramoto
  oscillators,'' in \emph{{E}uropean {C}ontrol {C}onference}, 2021, pp.
  1497--1502.

\bibitem{SA-DRM:2021a}
A.~Salova and R.~M. D'Souza, ``Cluster synchronization on hypergraphs,''
  \emph{arXiv preprint arXiv:2101.05464}, 2021.

\bibitem{SA-DRM:2021b}
------, ``Analyzing states beyond full synchronization on hypergraphs requires
  methods beyond projected networks,'' \emph{arXiv preprint arXiv:2107.13712},
  2021.

\bibitem{WW-WZ-TC:09}
W.~Wu, W.~Zhou, and T.~Chen, ``Cluster synchronization of linearly coupled
  complex networks under pinning control,'' \emph{IEEE {T}ransactions on
  {C}ircuits and {S}ystems}, vol.~56, no.~4, pp. 829--839, 2009.

\bibitem{GLV-FM-LV:2018}
L.~V. Gambuzza, M.~Frasca, and V.~Latora, ``Distributed control of
  synchronization of a group of network nodes,'' \emph{IEEE Transactions on
  Automatic Control}, vol.~64, no.~1, pp. 365--372, 2018.

\bibitem{DF-GR-MDB:17}
D.~Fiore, G.~Russo, and M.~di~Bernardo, ``Exploiting nodes symmetries to
  control synchronization and consensus patterns in multiagent systems,''
  \emph{Control Systems Letters}, vol.~1, no.~2, pp. 364--369, 2017.

\bibitem{TM-GB-DSB-FP:21-arxiv}
T.~Menara, G.~Baggio, D.~S. Bassett, and F.~Pasqualetti, ``Functional control
  of oscillator networks,'' \emph{arXiv:2012.04217}, 2021, submitted.

\bibitem{YQ-DSB-FP:22a}
Y.~Qin, D.~S. Bassett, and F.~Pasqualetti, ``Vibrational control of cluster
  synchronization: Connections with deep brain stimulation,'' in \emph{{IEEE}
  Conf.\ on Decision and Control}, Canc\'un, Mexico, Dec. 2022.

\bibitem{LT-CF-MI-DSB-FP:17}
L.~Tiberi, C.~Favaretto, M.~Innocenti, D.~S. Bassett, and F.~Pasqualetti,
  ``Synchronization patterns in networks of {K}uramoto oscillators: A geometric
  approach for analysis and control,'' in \emph{{IEEE} Conf.\ on Decision and
  Control}, Melbourne, Australia, Dec. 2017, pp. 481--486.

\bibitem{REB-JB-SMM:86a}
R.~E. Bellman, J.~Bentsman, and S.~M. Meerkov, ``Vibrational control of
  nonlinear systems: Vibrational stabilization,'' \emph{IEEE Transactions on
  Automatic Control}, vol.~31, no.~8, pp. 710--716, 1986.

\bibitem{SMM:80}
S.~M. Meerkov, ``Principle of vibrational control: Theory and applications,''
  \emph{IEEE Transactions on Automatic Control}, vol.~25, no.~4, 1980.

\bibitem{PRV-TM:80}
R.~Patel and M.~Toda, ``Quantitative measures of robustness for multivariable
  systems,'' in \emph{Joint Automation Control Conference}, no.~17, 1980,
  p.~35.

\bibitem{YR:85}
R.~Yedavalli, ``Improved measures of stability robustness for linear state
  space models,'' \emph{IEEE Transactions on Automatic Control}, vol.~30,
  no.~6, pp. 577--579, 1985.

\bibitem{haddad2011nonlinear}
W.~M. Haddad and V.~Chellaboina, \emph{Nonlinear {D}ynamical {S}ystems and
  {C}ontrol: {A} Lyapunov-{B}ased {A}pproach}.\hskip 1em plus 0.5em minus
  0.4em\relax Princeton, NJ, USA: Princeton University Press, 2011.

\bibitem{QY-KY-BDOA-CM:2021}
Y.~Qin, Y.~Kawano, B.~D. Anderson, and M.~Cao, ``Partial exponential stability
  analysis of slow-fast systems via periodic averaging,'' \emph{IEEE
  Transactions on Automatic Control}, 2021.

\bibitem{HKK:02-bis}
H.~K. Khalil, \emph{Nonlinear Systems}.\hskip 1em plus 0.5em minus 0.4em\relax
  Prentice Hall, 2002.

\bibitem{JBJ-GG:00}
J.~Bang-Jensen and G.~Gutin, \emph{Digraphs: Theory, Algorithms and
  Applications}, ser. Monographs in Mathematics.\hskip 1em plus 0.5em minus
  0.4em\relax Springer, 2000.

\bibitem{AMN-YQ-CAA-DSB-FP:22}
A.~M. Nobili, Y.~Qin, C.~A. Avizzano, D.~S. Bassett, and F.~Pasqualetti,
  ``Vibrational stabilization of complex network systems,'' in \emph{{A}merican
  {C}ontrol {C}onference}, San Diego, CA, May 2022.

\end{thebibliography}

\end{document}